\documentclass{jT}

\usepackage{amssymb,amsmath,latexsym,amsthm}%or any usual package

\newtheorem{theorem}{Theorem}[section]

\newtheorem{lemma}[theorem]{Lemma}

\newtheorem{proposition}[theorem]{Proposition}

\theoremstyle{definition}

\newtheorem*{example}{Example}

\begin{document}

\title[Weak Gibbs property]{Weak Gibbs property and systems of numeration}

% premier auteur

\author[\'Eric {\sc Olivier}]{{\sc \'Eric} OLIVIER}

\address{\'Eric {\sc Olivier}\\
Centre de Ressources Informatiques\\
Universit\'e de Provence\\
3, place Victor Hugo\\
13331 MARSEILLE Cedex 3, France}

\email{Eric.Olivier@up.univ-mrs.fr}

% deuxi\`eme auteur

\author[Alain {\sc Thomas}]{{\sc Alain} THOMAS}
\address{Alain {\sc Thomas}\\
Centre de Math\'ematiques et d'Informatique, LATP \'Equipe de th\'eorie des nombres\\
39, rue F. Joliot-Curie\\
13453 Marseille Cedex 13, France}

\email{thomas@cmi.univ-mrs.fr}

\maketitle

\begin{resume}

Nous \'etudions les propri\'et\'es d'autosimilarit\'e et la nature gibbsienne de certaines mesures definies sur l'espace produit
$\Omega_r:=\{0,1,\dots,r-1\}^{\mathbb N}$. Cet espace peut \^etre identifi\'e \`a l'intervalle $[0,1]$ au moyen de la num\'eration en
base $r$. Le dernier paragraphe concerne la convolution de Bernoulli en base $\beta={1+\sqrt5\over2}$, appel\'ee mesure de Erd\H os,
et son analogue en base $-\beta=-{1+\sqrt5\over2}$, que nous \'etudions au moyen d'un syst\`eme de numeration appropri\'e.

\end{resume}

\begin{abstr}

We study the selfsimilarity and the Gibbs properties of several measures defined on 
the product space 
$\Omega_r:=\{0,1,\dots,\break r-1\}^{\mathbb N}$. This space can be identified with the interval $[0,1]$ by means of the numeration in
base $r$. The last section is devoted to the Bernoulli convolution in base $\beta={1+\sqrt5\over2}$, called the Erd\H os measure,
and its analogue in base $-\beta=-{1+\sqrt5\over2}$, that we study by means of a suitable system of numeration.

\end{abstr}

\bigskip

\par\par

\bigskip

\section{Introduction}

One calls the Bernoulli convolution associated with the base $\beta>1$ and the parameter vector ${\bf p}=(p_0,\dots,
p_{s-1})$, the infinite product  of the Dirac measures 
$p_0\delta_{0\over \beta^n}+\dots+p_{s-1}\delta_{s-1\over \beta^n}$ for $n\ge1$ (see \cite{Erd,SV,OST2,OST3}). In other words, it is the distribution measure
of the random variable defined~by
$$
X(\omega)=\sum_{n\ge1}{\omega_n\over \beta^n},
$$
when $\omega=(\omega_n)_{n\in\mathbb N}$ has a Bernoulli distribution such that, for any $n\in{\mathbb N}$,$$P(\omega_n=0)=p_0,\dots,P(\omega_n=s-1)=p_{s-1}.$$
The Bernoulli convolution associated with $\beta$ and ${\bf p}$ is the unique measure $\mu$ with bounded support that satisfies the
self-similarity relation (\cite{PS}):
$$
\mu=\sum_{i=0}^{s-1}p_i\cdot\mu\circ
S_i^{-1},
$$
where the affine contractions $S_i:{\mathbb R}\to{\mathbb R}$ are defined
 by $S_i(x):={x+i\over\beta}$.
The measure $\mu$ is purely singular with respect to the Lebesgue measure when ${\bf p}$ is uniform and $\beta$ a
Pisot number, that is, the conjuguates of $\beta$ have modulus less than $1$.
The problem to know if $\mu$ has the weak Gibbs property in the sense of Yuri \cite{Yur} is not simple; it is solved in case $\beta$ is
a multinacci number (\cite{FO,OST3}), but more complicated for other Pisot numbers of degree at least~3
(for instance in \cite[Example 2.4]{OST3}, computing the values of the Bernoulli convolution in case $\beta^3=3\beta^2-1$
requires matrices of order $8$).

Section 2 recalls the definition of the weak Gibbs property, and its link with the notions of Bernoulli or Markov measure.

Section 3 is devoted to some results of Mukherjea, Nakassis and Ratti about products of i. i. d. random stochastic matrices, 
that we present in a slighty different way 
(Proposition \ref{infiniteproduct}). 
They have computed the density of the limit distribution, in case this distribution is the 
Bernoulli convolution in base $\beta=\root m\of r$ with parameters $p_0=\dots=p_{r-1} ={1\over r}$.

The framework is different in the sections 5 to 7; we define a measure on $\Omega_r:=\{0,1,\dots,r-1\}^{\mathbb N}$ 
by giving its values on the cylinders of $\Omega_r$, under the form of products of $2\times2$ matrices and vectors. Theorem \ref{theorem} gives 
a condition for such a measure,
to be related to a Bernoulli convolution, via the representation of the reals in the integral base $r$.
Establishing the weak Gibbs property requires the convergence of the involved product of matrices and vectors in the projective space of dimension $2$.
It is proved in \cite{FO} that the uniform Bernoulli convolution in base $\beta={1+\sqrt5\over2}$ is weak Gibbs; here, Section 7 give analogue result 
in the base $-\beta=-{1+\sqrt5\over2}$.

\section{Weak Gibbs measures}

One says that the probability measure $\mu$ on the product space 
$\Omega_r=\{0,1,\dots,r-~1\}^{\mathbb N}$ has the
weak Gibbs property if there exists a map $\phi:~\Omega_r\to{\mathbb R}$, continuous for the product topology on $\Omega_r$,
such that
\begin{equation}\label{wG}
\lim_{n\to\infty}\left({\mu[\omega_1\dots\omega_n]\over
e^{\phi(\omega)}e^{\phi(\sigma\omega)}\dots e^{\phi(\sigma^{n-1}\omega)}}\right)^{1/n}=1\quad\hbox{uniformly  on
}\omega\in\Omega_r
\end{equation}
(where $\sigma$ is the shift on $\Omega_r$~, and $[\omega_1\dots\omega_n]$ is the cylinder of order $n$ around $\omega$ that is,
the set of the
$\omega'\in\Omega_r$ such that $\omega'_i=\omega_i$ for $1\le i\le n$). 
If (\ref{wG})~holds, $\phi$ is called a potential of $\mu$.

Equivalently, $\mu$ has the
weak Gibbs property if and only if the measure of any cylinder $[\omega_1\dots\omega_n]$ can be approached by a product in the
following way: there exists a continuous map \hbox{$\varphi:\Omega_r\to]0,+\infty[$} such that
$$
\begin{array}{c}\forall\varepsilon>0\enskip\exists N\in{\mathbb N}\enskip\forall n\ge N\enskip
\forall\omega\in\Omega_r\\
(\varphi(\omega)-\varepsilon)\dots(\varphi(\sigma^{n-1}\omega)-\varepsilon)\le
\mu[\omega_1\dots\omega_n]\le(\varphi(\omega)+\varepsilon)\dots(\varphi(\sigma^{n-1}\omega)+\varepsilon).
\end{array}
$$

In case $\mu$ is $\sigma$-invariant, the following theorem gives an equivalent definition (see \cite{Kea}, \cite{Wal},
\cite{P-P-W}), which involves the map $\phi_\mu$ defined as follows:
\begin{equation}\label{potential}
\phi_\mu(\omega):=\lim_{n\to\infty}\log{\mu[\omega_1\dots\omega_n]\over\mu[\omega_2\dots\omega_n]}
\end{equation}
at each point $\omega\in\Omega_r$ such that this limit exists.

\begin{theorem}

Let $\mu$ be a 
$\sigma$-invariant probability measure on $\Omega_r$~, and \hbox{$\phi:\Omega_r\to{\mathbb R}$} a continuous map.
The following assertions are equivalent:

(i) $\mu$ is a weak Gibbs measure of potential $\phi$ and, for any $\omega\in\Omega_r$~,
$\displaystyle\sum_{a=0}^{r-1}e^{\phi(a\omega)}=1$;

(ii) $\phi_\mu(\omega)$ exists for any $\omega\in\Omega_r$~, and
$\phi_\mu=\phi$;

(iii) $\mu$ has entropy $h_\mu=-\mu(\phi)$ and, for any $\omega\in\Omega_r$~,
$\displaystyle\sum_{a=0}^{r-1}e^{\phi(a\omega)}=1$.

\end{theorem}

This theorem can be used to prove that a $\sigma$-invariant probability measure has the weak Gibbs property, by using the implication
$(ii)\Rightarrow(i)$. Now for any probability measure $\mu$ on $\Omega_r$, not necessarily $\sigma$-invariant,
the following implication is straightforward (see \cite{OST3}):
\begin{proposition}\label{pot=>weak}
If $\phi_\mu$ is defined and continuous on $\Omega_r$~, then $\mu$ is a weak Gibbs measure of potential
$\phi_\mu$.
\end{proposition}

The two following examples show that the Bernoulli and the Markovian measures are weak Gibbs.
The third is a counterexample: the potential of the weak Gibbs measure $\mu_3$ is not $\phi_{\mu_3}$.

\begin{example}

If $\mu_1$ is a Bernoulli measure
with support
$\Omega_r$~, then $\phi_{\mu_1}$ is the continuous map such that 
$\phi_{\mu_1}(\omega)=\log\mu_1[\omega_1]$ for any $\omega\in\Omega_r$~.

\end{example}

\begin{example}If $\mu_2$ is a Markov measure with support $\Omega_r$~, then $\phi_{\mu_2}$ is the continuous map such that
$\displaystyle\phi_{\mu_2}(\omega)=\log {\mu_2[\omega_1\omega_2]\over\mu_2[\omega_2]}$ for any
$\omega\in\Omega_r$~.
%(notice that the logarithms of the passage probabilities
%provide another potential for $\mu_2$)

\end{example}

\begin{example}

(see \cite{OST2}) Let the probability measure $\mu_3$ be defined on
$\Omega_r$ by
$$
\mu_3[\omega_1\dots\omega_n]:={1\over2\cdot(2r)^n}\left(\begin{array}{cc}1&1\end{array}\right)
\left(\begin{array}{cc}\omega'_1&0\cr1&1\end{array}\right)\dots
\left(\begin{array}{cc}\omega'_n&0\cr1&1\end{array}\right)\left(\begin{array}{cc}1\cr1\end{array}\right)
$$
with $\omega'_i=1+{2\omega_i\over r-1}$. Then $\mu_3$ is weak Gibbs of potential $\phi:\omega\mapsto\log{\omega'_1\over2r}$,
although $\phi_{\mu_3}$ is discontinuous at any point $\omega$ such that the series\hfil\break$S_\omega:=\sum_{n\ge1}{1\over\omega'_1\dots\omega'_n}$ converges:
$$
\phi_{\mu_3}(\omega)=\left\{\begin{array}{ll}\log{1\over2r}+\log(1+{1\over S_\omega})&\hbox{if }
S_\omega<\infty\\
\log{1\over2r}&\hbox{if }
S_\omega=\infty.\end{array}\right.
$$

\end{example}

The notion of weak Gibbs measure generalize the one of Gibbs measure (see for instance \cite{Bow1}). Let us generalize in the
same way the notion of quasi-Bernoulli measure (see
\cite{BMP}, \cite{Heu}), and say that
$\mu$ is weakly quasi-Bernoulli~if it satisfies the following condition:
\begin{equation}\label{wasi-Bernoulli}
\lim_{n\to\infty}\left({\mu[\omega_1\dots\omega_n]\over
\mu[\omega_1\dots\omega_i]\mu[\omega_{i+1}\dots\omega_n]}\right)^{1/n}=1\quad\hbox{uniformly  on
}\omega\in\Omega_r\hbox{ and }i\in\{1,\dots,n\}.
\end{equation}

Then one has the following

\begin{proposition}\label{wasi-proposition}If a probability measure on $\Omega_r$ has the weak Gibbs property, it satisfies (\ref{wasi-Bernoulli}).
\end{proposition}
This proposition is straightforward, but can be used to prove that a probability measure do not have the weak Gibbs property:

\begin{example}  Let $\mu_4$ be defined on
$\Omega_2$ by
$$
\mu_4[\omega_1\dots\omega_n]:={1\over2}\left(\begin{array}{cc}1&1\end{array}\right)M_{\omega_1}\dots
M_{\omega_n}\left(\begin{array}{cc}1\cr1\end{array}\right)
$$where $M_0={1\over4}\left(\begin{array}{cc}0&0\cr1&1\end{array}\right)$ and $M_1={1\over4}\left(\begin{array}{cc}4&0\cr1&1\end{array}\right)$.
It is not weak Gibbs because $\displaystyle\left({\mu_4[1^n0^n]\over
\mu_4[1^n]\mu_4[0^n]}\right)^{1/n}$ do not converge to 1.
\end{example}
One can ask if the converse of Proposition \ref{wasi-proposition} true, or if the condition (\ref{wasi-Bernoulli}) imply
that $\mu$ is weak Bernoulli in the sense of Bowen \cite{Bow}.

\section{Products of stochastic matrices}

We consider a finite set of stochastic $2\times2$ matrices, let $M_k=\left(\begin{array}{cc}x_k&1-x_k\cr y_k&1-y_k\end{array}\right)$ for $k=0,1,\dots,r-1$,
where $x_k,y_k\in[0,1]$. We suppose the $M_k$ are different from $\left(\begin{array}{cc}1&0\cr0&1\end{array}\right)$ and
$\left(\begin{array}{cc}0&1\cr 1&0\end{array}\right)$.

A. Mukherjea and al. have studied in \cite{MNR} and \cite{MN} the distribution of the random matrice 
\hbox{$\Omega_r\ni\omega\mapsto M_{\omega_1}\dots M_{\omega_n}$} when the distribution of $\omega$ is Bernoulli 
with positive parameters $p_0,\dots p_{r-1}$. This distribution converges
when $n\to\infty$, though the matrix $M_{\omega_1}\dots M_{\omega_n}$ itself
do not converge (that is, its entries are -- in much cases -- divergent sequences). But we shall prove the
convergence of the matrix $M_{\omega_n}\dots M_{\omega_1}$ (which has of course the same distribution as 
$M_{\omega_1}\dots M_{\omega_n}$ when the distribution of $\omega$ is Bernoulli).
\begin{proposition}\label{infiniteproduct}The product matrix $P_n^\omega:=M_{\omega_n}\dots M_{\omega_1}$ 
converges uniformly on $\omega\in\Omega_r$ to the matrix 
$\left(\begin{array}{cc}x^\omega&1-x^\omega\cr x^\omega&1-x^\omega\end{array}\right)$,
where $x^\omega:=\sum_{i=1}^{\infty}y_{\omega_i}\det P_{i-1}^\omega$
and -- by convention -- $P_0^\omega$ is the unit-matrix.
\end{proposition}
\begin{proof} Setting $x_n^\omega:=y_n^\omega+\det P_n^\omega$ 
with $y_n^\omega:=\sum_{i=1}^ny_{\omega_i}\det P_{i-1}^\omega$ one check easily by induction that
$$
P_n^\omega=\left(\begin{array}{cc}x_n^\omega&1-x_n^\omega\cr y_n^\omega&1-y_n^\omega\end{array}\right).
$$
The uniform convergence of the sequences $x_n^\omega$ and $y_n^\omega$ is due to the fact that 
each matrix $M_k$ has -- from the hypotheses -- a determinant less than 1 in absolute value.

\end{proof}

\begin{theorem} (\cite[Section 2]{MNR}) The distribution of $\omega\mapsto x^\omega$ is\hfil\break
\_ discrete if at least one of the matrices $M_k$ is non invertible;\hfil\break
\_ singular continuous if the product $\left({|\det M_0|\over p_0}\right)^{p_0}\dots\left({|\det
M_{r-1}|\over p_{r-1}}\right)^{p_{r-1}}$ belongs to $]0,1]$ and at least one of its factors is different from $1$.
\end{theorem}

\underline{Selfsimilarity relation:} The random variable $\omega\mapsto x^\omega$ takes its
values in $[0,1]$ because $\left(\begin{array}{cc}x^\omega&1-x^\omega\cr x^\omega&1-x^\omega\end{array}\right)$ is the limit of nonnegative matrices.
Let $\lambda$ be the probability distribution of $\omega\mapsto x^\omega$. If all the matrices $M_k$ are invertible,
then $\lambda$ is selfsimilar in the sense that, for any borelian $B\subset[0,1]$,

$$\lambda(B)=\sum_{k=0}^{r-1}p_k\ \lambda\left({B-y_k\over x_k-y_k}\right)$$

(see \cite[equation~(2.6)]{MNR} for the proof).

Let us represent, for instance in the case $r=2$ with $(x_0-y_0)(x_1-y_1)<0$, 
the two maps $\displaystyle x\mapsto{x-y_k\over x_k-y_k}$ involved in the selfsimilarity
relation:

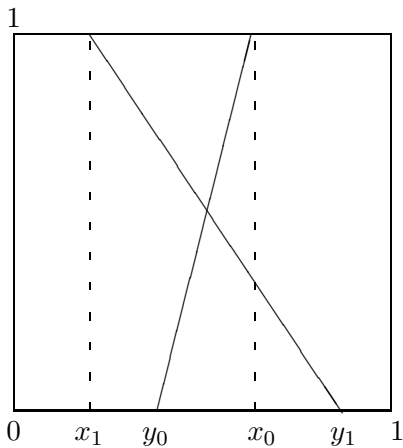
\setlength{\unitlength}{1cm}
\begin{picture}(6.5,6)(5,5)
\put(5,5){\framebox(5,5)}
\put(6.9,5){\line(1,4){1.25}}
\put(6,10){\line(2,-3){3.36}}
\put(4.9,10.1){$1$}
\put(4.9,4.6){$0$}
\put(8.1,4.6){$x_0$}
\put(5.8,4.6){$x_1$}
\put(8.2,5){\line(0,1){0.1}}
\put(6,5){\line(0,1){0.1}}
\put(8.2,5.4){\line(0,1){0.1}}
\put(6,5.4){\line(0,1){0.1}}
\put(8.2,5.8){\line(0,1){0.1}}
\put(6,5.8){\line(0,1){0.1}}
\put(8.2,6.2){\line(0,1){0.1}}
\put(6,6.2){\line(0,1){0.1}}
\put(8.2,6.6){\line(0,1){0.1}}
\put(6,6.6){\line(0,1){0.1}}
\put(8.2,7){\line(0,1){0.1}}
\put(6,7){\line(0,1){0.1}}
\put(8.2,7.4){\line(0,1){0.1}}
\put(6,7.4){\line(0,1){0.1}}
\put(8.2,7.8){\line(0,1){0.1}}
\put(6,7.8){\line(0,1){0.1}}
\put(8.2,8.2){\line(0,1){0.1}}
\put(6,8.2){\line(0,1){0.1}}
\put(8.2,8.6){\line(0,1){0.1}}
\put(6,8.6){\line(0,1){0.1}}
\put(8.2,9){\line(0,1){0.1}}
\put(6,9){\line(0,1){0.1}}
\put(8.2,9.4){\line(0,1){0.1}}
\put(6,9.4){\line(0,1){0.1}}
\put(8.2,9.8){\line(0,1){0.1}}
\put(6,9.8){\line(0,1){0.1}}
\put(6.7,4.6){$y_0$}
\put(9.2,4.6){$y_1$}
\put(10,4.6){$1$}
\end{picture}

\bigskip

\begin{example} The probability distribution $\lambda$ of $\omega\mapsto x^\omega$ is related to the numeration in a given
base $\beta>1$ if we suppose that $x_k=y_k+{1\over\beta}$ and that $y_0,\dots,y_{r-1}$ are in arithmetic progression.
Since we want that $x_k$ and $y_k$ belong to $[0,1]$, the good choice is
$$
y_k={k\over r-1}\left(1-{1\over\beta}\right)\quad\hbox{for }k=0,\dots,r-1.
$$
Then $x_\omega={\beta-1\over r-1}\sum_{n\ge1}{\omega_n\over
\beta^n}$ and $\lambda\left({\beta-1\over r-1}\cdot\right)$ is the convolution of the measures 
$p_0\delta_{0\over \beta^n}+\dots+
p_{r-1}\delta_{r-1\over \beta^n}$ for $n=1,2,\dots$.

In case $\beta=\root m\of r$ with $m\in{\mathbb N}$, if the distribution is uniform
($p_0=\dots=p_{r-1} ={1\over r}$), it is proved in
\cite[Proposition 1]{MNR} that the density of the (absolutely continuous) distribution of
$\omega\mapsto x^\omega$ is a piece-wise polynomial of degree at most~$m$.\end{example}

\section{Uniform convergence (in direction) of the sequence of vectors $n\mapsto M_{\omega_1}\dots M_{\omega_n}V$}

In this section ${\mathcal M}=\{M_0,\dots,M_{r-1}\}$ is a finite set of $2\times2$ matrices, 
where each matrix
$M=\left(\begin{array}{cc}a&b\cr c&d\end{array}\right)\in{\mathcal M}$ has nonnegative entries and each of the
columns $\left(\begin{array}{c}a\cr c\end{array}\right)$ and $\left(\begin{array}{c}b\cr d\end{array}\right)$ is distinct 
from~$\left(\begin{array}{c}0\cr0\end{array}\right)$.
One denotes by ${\mathcal M}^2$ the set of matrices $MM'$ for $M$ and $M'$ in ${\mathcal M}$, and\hfil\break
${\mathcal M}_1$ the set of matrices $M\in{\mathcal M}$ with
$a=0$;\hfil\break${\mathcal M}_2$ the set of matrices in $M\in{\mathcal M}^2$ with $b=0$;\hfil\break${\mathcal M}_3$ the set of matrices in $M\in{\mathcal M}^2$ with
$c=0$;\hfil\break 
${\mathcal M}_4$ the set of matrices in $M\in{\mathcal M}$ with $d=0$.

\begin{proposition}(\cite[theorem A]{OST2}) Let $V=\left(\begin{array}{cc}v_1\cr v_2\end{array}\right)$ be
a column matrix
with positive entries. The sequence $\displaystyle n\mapsto{M_{\omega_1}\dots M_{\omega_n}V
\over\|M_{\omega_1}\dots M_{\omega_n}V\|}$ converges uniformly on $\omega\in\Omega_r$ if and only if at least one of the
following conditions holds:

$(i)\ \not\exists  M\in{\mathcal M}_2$ such that $a>d$ and $\not\exists  M\in{\mathcal M}_3$ such that $a<d$ and 
${\mathcal M}_2\cap{\mathcal M}_3=\emptyset$

$(ii)\ \not\exists M\in{\mathcal M}_2$ such that $a\le d$, and $\not\exists  M\in{\mathcal M}_3$ such that $a\ge d$

$(iii)\ \not\exists  M\in{\mathcal M}_2$ such that $a\le d$ and $\not\exists  M\in{\mathcal M}_3$ such that $a<d$ and ${\mathcal
M}_1=\emptyset$

$(iv)\ \not\exists  M\in{\mathcal M}_2$ such that $a>d$ and $\not\exists  M\in{\mathcal M}_3$ such that $a\ge d$ and ${\mathcal
M}_4=\emptyset$

$(v)\ V$ is an eigenvector of all the matrices in $\mathcal M$.
\end{proposition}

\section{Application to the measures defined by products of
matrices}

Let ${\mathcal M}=\{M_0,\dots,M_{r-1}\}$ be a finite set of $2\times2$ matrices whose columns are distinct
from~$\left(\begin{array}{cc}0\cr0\end{array}\right)$, and let $L$ (resp. $V$) be a positive row matrix (resp., a positive column matrix). If $V$ is an eigenvector of
$\sum_iM_i$ for the eigenvalue $1$, one can define some measure $\eta$ on $\Omega_r$ by setting
$$
\eta[\omega_1\dots\omega_n]=LM_{\omega_1}\dots M_{\omega_n}V.
$$

\begin{proposition}\label{siconditions}The map $\phi_\eta$ defined in (\ref{potential}), exists and is continuous if and only if ${\mathcal M}$ satisfies at least one of
the above conditions $(i)$,
\dots, $(v)$, or the following:

$(vi)\ L$  is an eigenvector of all the matrices in ${\mathcal M}$.
\end{proposition}

\begin{proof} The map $\phi_\eta$ is related to the map $\displaystyle\psi_{\mathcal M}:\omega\mapsto
\lim_{n\to\infty}{M_{\omega_1}\dots
M_{\omega_n}V
\over\|M_{\omega_1}\dots M_{\omega_n}V\|}$~. Indeed
$$
\phi_\eta(\omega)={LM_{\omega_1}\psi_{\mathcal M}\circ\sigma(\omega)\over L\psi_{\mathcal M}\circ\sigma(\omega)}
$$
for any $\omega\in\Omega_r$ such that
$\psi_{\mathcal M}\circ\sigma(\omega)$ exists. Moreover if $(vi)$ does not hold, the domains of
definition of $\phi_\eta$ and $\psi_{\mathcal M}\circ\sigma$ are the same.

\end{proof}

Now this proposition gives a sufficient condition for 
$\eta$ to have the weak Gibbs property (by using Proposition \ref{pot=>weak}). 
This condition is of course not necessary (see Example 1.5).

\section{Measures associated with the numeration in integral base $r$.}

Let the map $\displaystyle X_{q,r}:\Omega_q\mapsto\left[0,{q-1\over r-1}\right]$ be defined by
$$
X_{q,r}(\omega)=\sum_{n\ge1}{\omega_n\over r^n}.
$$
In particular $X_{r,r}$ is one-to-one except on a countable set because, if $\omega$ is not eventually $r-1$, 
the real $X_{r,r}(\omega)$ has expansion $\omega$ in base $r$. In the present section we identify 
the set of sequences $\Omega_r$ with the interval $[0,1]$, by means the map $X_{r,r}$~.

The following theorem gives a condition for a measure defined by products of $2\times2$ matrices,
to be related to some Bernoulli convolution in base~$r$:

\begin{theorem}\label{theorem}(\cite{OST2}, Theorem 4.25) Let $\nu$ be a $\sigma$-invariant probability measure on $\Omega_r$~; the
following assertions are equivalent:

$(i)$ there exists a nonnegative row matrix $L$, a column matrix $V$ and some square matrices
$M_0,\dots,M_{r-1}$ such that
$$
\forall\omega\in\Omega_r,\ \forall n\in{\mathbb N},\quad
\nu[\omega_1\dots\omega_n]=LM_{\omega_1}\dots M_{\omega_n}V,
$$
where the matrices $M_k=\left(\begin{array}{cc}a_k&b_k\cr c_k&d_k\end{array}\right)$ satisfy the conditions
$$
b_0=0\hbox{ and }\left(\begin{array}{c}a_k\cr c_k\end{array}\right)=\left\{\begin{array}{ll}
\left(\begin{array}{c}b_{k+1}\cr d_{k+1}\end{array}\right)&\hbox{if }0\le k\le r-2\\
\left(\begin{array}{c}d_0\cr b_0\end{array}\right)&\hbox{if }k=r-1
\end{array}\right.
$$

$(ii)$ there exists a nonnegative parameter vector ${\bf p}=(p_0,\dots,
p_{2r-2})$ such that $\nu$ is the probability distribution $\nu_{\bf p}$ of the fractional part of 
$X_{2r-1,r}(\omega)$,
when $\omega\in\Omega_{2r-1}$ has a Bernoulli distribution with 
parameter $\bf p$.
\end{theorem}
The relations between the matrices $M_k$ and the parameter $\bf p$ are 
$$
p_0=a_0,\dots,
p_{r-1}=a_{r-1},p_r=c_0,\dots, p_{2r-2}=c_{r-2}
$$
and thus $\nu_{\bf p}$ is weak Gibbs from Proposition \ref{siconditions} in certain cases, for instance if the $p_k$ are positive.

\underline{Selfimilarity relation} Let $\mu_{\bf p}$ and $\nu_{\bf p}$ be
the probability distributions of $X_{2r-1,r}$ and the fractionnal part of 
$X_{2r-1,r}$, respectively. Their respective supports are $[0,2]$ and $[0,1]$ and,
for any borelian $B\subset[0,1]$,
$$
\nu_{\bf p}(B)=
\mu_{\bf p}(B)+
\mu_{\bf p}(B+1).
$$ 
Theorem \ref{theorem} is a consequence of the selfsimilarity relation
\begin{equation}\label{auto}
\mu_{\bf p}(B)=\sum_{k=0}^{2(r-1)}p_k\ \mu_{\bf p}(rB-k)
\end{equation}
which allows to compute the column matrix $\left(\begin{array}{c}\mu_{\bf p}(B)\cr
\mu_{\bf p}(B+1)\end{array}\right)$.

The measure $\nu_{\bf p}$ has support $[0,1]$, while the measure $\nu^*_{\bf p}$ defined 
for any borelian $B\subset{\mathbb R}$, by
$$
\nu^*_{\bf p}(B)=\mu_{\bf p}(B)+\mu_{\bf p}(B+1)
$$
has support $[-1,2]$, and coincide with $\nu_{\bf p}$ on $[0,1]$. 
The selfsimilarity relation for $\nu^*_{\bf p}$ can be deduced from (\ref{auto}):
\begin{equation}\label{autoetoile}
\nu^*_{\bf p}(B)=\sum_{k=-(r-1)}^{2(r-1)}p^*_k\ \nu^*_{\bf p} (rB-k)
\end{equation}
where $p^*_k=\sum_{j\ge0}(p_{k+2j}-
p_{k+2j+1}+p_{k+2j+b}-
p_{k+2j+b+1})$.

Both measures $\mu_{\bf p}$ and $\nu^*_{\bf p}$ are
Bernoulli convolutions: they are -- respectively -- the infinite product  of the measures 
$p_0\delta_{0\over r^n}+\dots+
p_{2r-2}\delta_{2r-2\over r^n}$ and the one of the measures 
$p^*_{-r+1}\delta_{{-r+1}\over r^n}+\dots+
p^*_{2r-2}\delta_{2r-2\over r^n}$~, for $n\ge1$.

We represent below the maps $x\mapsto rx-k$ involved in (\ref{auto}) and (\ref{autoetoile}), in the case $r=4$:

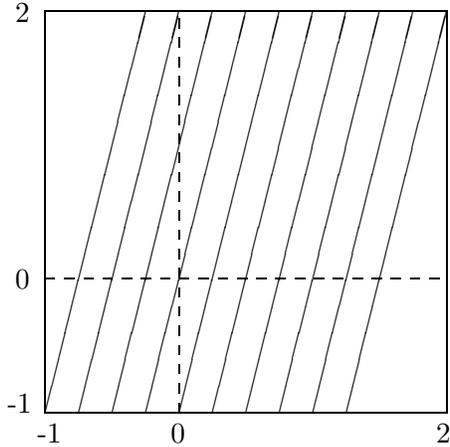
\setlength{\unitlength}{1cm}
\begin{picture}(6.5,6.5)(5,5)
\put(5,5){\framebox(5.33,5.33)}
\put(5,5){\line(1,4){1.333333}}
\put(5.444444,5){\line(1,4){1.333333}}
\put(5.888888,5){\line(1,4){1.333333}}
\put(6.333333,5){\line(1,4){1.333333}}
\put(6.777777,5){\line(1,4){1.333333}}
\put(7.222222,5){\line(1,4){1.333333}}
\put(7.666666,5){\line(1,4){1.333333}}
\put(8.111111,5){\line(1,4){1.333333}}
\put(8.555555,5){\line(1,4){1.333333}}
\put(9,5){\line(1,4){1.333333}}
\put(6.777777,5){\line(0,1){0.127}}
\put(6.777777,5.255){\line(0,1){0.127}}
\put(6.777777,5.511){\line(0,1){0.127}}
\put(6.777777,5.766){\line(0,1){0.127}}
\put(6.777777,6.022){\line(0,1){0.127}}
\put(6.777777,6.277){\line(0,1){0.127}}
\put(6.777777,6.533){\line(0,1){0.127}}
\put(6.777777,6.788){\line(0,1){0.127}}
\put(6.777777,7.044){\line(0,1){0.127}}
\put(6.777777,7.299){\line(0,1){0.127}}
\put(6.777777,7.555){\line(0,1){0.127}}
\put(6.777777,7.811){\line(0,1){0.127}}
\put(6.777777,8.066){\line(0,1){0.127}}
\put(6.777777,8.322){\line(0,1){0.127}}
\put(6.777777,8.577){\line(0,1){0.127}}
\put(6.777777,8.833){\line(0,1){0.127}}
\put(6.777777,9.088){\line(0,1){0.127}}
\put(6.777777,9.344){\line(0,1){0.127}}
\put(6.777777,9.599){\line(0,1){0.127}}
\put(6.777777,9.855){\line(0,1){0.127}}
\put(6.777777,10.111){\line(0,1){0.127}}

\put(5,6.777777){\line(1,0){0.127}}
\put(5.255,6.777777){\line(1,0){0.127}}
\put(5.511,6.777777){\line(1,0){0.127}}
\put(5.766,6.777777){\line(1,0){0.127}}
\put(6.022,6.777777){\line(1,0){0.127}}
\put(6.277,6.777777){\line(1,0){0.127}}
\put(6.533,6.777777){\line(1,0){0.127}}
\put(6.788,6.777777){\line(1,0){0.127}}
\put(7.044,6.777777){\line(1,0){0.127}}
\put(7.299,6.777777){\line(1,0){0.127}}
\put(7.555,6.777777){\line(1,0){0.127}}
\put(7.811,6.777777){\line(1,0){0.127}}
\put(8.066,6.777777){\line(1,0){0.127}}
\put(8.322,6.777777){\line(1,0){0.127}}
\put(8.577,6.777777){\line(1,0){0.127}}
\put(8.833,6.777777){\line(1,0){0.127}}
\put(9.088,6.777777){\line(1,0){0.127}}
\put(9.344,6.777777){\line(1,0){0.127}}
\put(9.599,6.777777){\line(1,0){0.127}}
\put(9.855,6.777777){\line(1,0){0.127}}
\put(10.11,6.777777){\line(1,0){0.127}}

\put(4.9,4.6){-$1$}
\put(6.67,4.6){$0$}
\put(10.2,4.6){$2$}
\put(4.5,5){-$1$}
\put(4.6,6.65){$0$}
\put(4.6,10.2){$2$}
\end{picture}

\section{The bases $\beta={1+\sqrt5\over2}$ and $-\beta=-{1+\sqrt5\over2}$}

We consider in this section the measures $\mu$ and $\mu_\star$
which are respectively the distributions of the random variables $X$ and $Y$, defined by
$$
X(\omega)=\sum_{n\ge1}{\omega_n\over \beta^{n+1}}\quad\hbox{and}\quad
Y(\omega)={1\over\beta}-\sum_{n\ge1}{\omega_n\over
(-\beta)^{n+1}},
$$
when the distribution of $\omega\in\Omega_2$ is Bernoulli with positive parameter vector ${\bf p}=(p,q)$.
We use consecutively two systems of numeration
(see for instance \cite{P}, \cite{OST} and \cite{DS}): any real $x\in[0,1[$ can be represented in an unique way on the form
$$
x=\sum_{n\ge1}{\varepsilon_n\over \beta^n}\quad\hbox{(Parry expansion)}\quad\hbox{and}\quad
x={1\over\beta}-\sum_{n\ge1}{\alpha_n\over
(-\beta)^{n+1}}\ ,
$$
where $(\varepsilon_n)_{n\ge1}=:\varepsilon(x)$ and
$(\alpha_n)_{n\ge1}=:\alpha(x)$ are two sequences with terms in $\{0,1\}$, 
without two consecutive terms $1$, such that
$\sigma^n\varepsilon(x)$ and
$\sigma^{2n+1}\alpha(x)$ differ from the periodic sequence $1010\dots$ for any $n\ge0$.
For any word $w=\omega_1\dots\omega_n$ on the alphabet $\{0,1\}$ and without factor $11$,
we denote
$$
\begin{array}{rcl}
[\![w]\!]&:=&\left\{x\in[0,1[\;;\;\varepsilon(x)\in[\omega_1,\dots,\omega_n]\right\}
\cr[\![w]\!]_\star&:=&\left\{x\in[0,1[\;;\;\alpha(x)\in[\omega_1,\dots,\omega_n]\right\}.
\end{array}
$$
In case $\omega_n=0$ we may compute $\mu[\![w]\!]$ and
$\mu_\star[\![w]\!]_\star$ by the following formulas:
\begin{equation}\label{mu}
\left(\begin{array}{c}\mu({1\over\beta}[\![w]\!])\cr\mu({1\over\beta}+{1\over\beta}[\![w]\!])\cr
\mu({1\over\beta^2}+{1\over\beta}[\![w]\!])\end{array}\right)=M_{\omega_1}\dots M_{\omega_n}\left(\begin{array}{c}{p\over
p+q^2}\cr{q^2\over
p+q^2}\cr{q\over
p+q^2}\end{array}\right)
\end{equation}
\begin{equation}\label{mustar}
\left(\begin{array}{c}\mu_\star([\![w]\!]_\star)\cr
\mu_\star(-{1\over\beta}+[\![w]\!]_\star)\cr
\mu_\star({1\over\beta^2}+[\![w]\!]_\star)\end{array}\right)=A_{\omega_1}\dots A_{\omega_n}\left(\begin{array}{c}1
\cr{q\over
1+q}\cr{1\over
1+q}\end{array}\right)
\end{equation}
where
$$
M_0=\left(\begin{array}{ccc}p&0&0\cr0&0&q\cr q&p&0\end{array}\right)\ 
M_1=\left(\begin{array}{ccc}q&p&0\cr0&0&0\cr0&q&0\end{array}\right)\ 
A_0=\left(\begin{array}{ccc}p&q&0\cr0&0&q\cr0&p&0\end{array}\right)\ 
A_1=\left(\begin{array}{ccc}q&0&0\cr0&0&0\cr p&q&0\end{array}\right).
$$
The formula (\ref{mu}) -- and its extension to the multinacci bases -- is proved in \cite{OST3}.
Let us sketch the proof of (\ref{mustar}), which is equivalent to the following (assuming again that the word $w$ do not have 
two consecutive letters $1$ and ends by the letter $0$):
\begin{equation}\label{mustarstep}
\left(\begin{array}{c}\mu_\star([\![w]\!]_\star)\cr
\mu_\star(-{1\over\beta}+[\![w]\!]_\star)\cr
\mu_\star({1\over\beta^2}+[\![w]\!]_\star)\end{array}\right)=A_{\omega_1}\left(\begin{array}{c}\mu_\star([\![w']\!]_\star)\cr
\mu_\star(-{1\over\beta}+[\![w']\!]_\star)\cr
\mu_\star({1\over\beta^2}+[\![w']\!]_\star)\end{array}\right)\quad\hbox{and}\quad\left(\begin{array}{c}\mu([0,1[)
\cr\mu(-{1\over\beta}+[0,1[)\cr\mu({1\over\beta^2}+[0,1[)\end{array}\right)=\left(\begin{array}{c}1
\cr{q\over
1+q}\cr{1\over
1+q}\end{array}\right)
\end{equation}
where $w'=\omega_2\dots\omega_n$ for $n\ge2$ and, by convention, if $n=1$ the word
$w'$ is empty and $[\![w']\!]_\star=[0,1[$. We first compute $\mu_\star([\![w]\!]_\star)$: it is the probability of the event
$Y(\xi)\in[\![w]\!]_\star$. This event is equivalent to 
$Y(\sigma\xi)\in{\omega_1-\xi_1\over\beta}+[\![w']\!]_\star$ hence

$\bullet$ in case $\omega_1=0$, it is also equivalent to
$$\left\{\begin{array}{l}\xi_1=0\\Y(\sigma\xi)\in[\![w']\!]_\star\end{array}\right.
\quad\hbox{or}\quad\left\{\begin{array}{l}\xi_1=1\\Y(\sigma\xi)\in-{1\over\beta}+[\![w']\!]_\star\end{array}\right.$$
and this explain why the first row in $A_0$ is $\left(\begin{array}{ccc}p&q&0\end{array}\right)$;

$\bullet$ in case $\omega_1=1$ we have necessarily $n\ge2$ and $\omega_2=0$, and the event
$Y(\sigma\xi)\in{1-\xi_1\over\beta}+[\![w']\!]_\star$ can occur only if $\xi_1=1$ and
$Y(\sigma\xi)\in[\![w']\!]_\star$; so the first row in $A_1$ is $\left(\begin{array}{ccc}q&0&0\end{array}\right)$. 

$\bullet$ We compute in the same way $\mu_\star(-{1\over\beta}+[\![w]\!]_\star)$ and
$\mu_\star({1\over\beta^2}+[\![w]\!]_\star)$ and we conclude that the first equality in
(\ref{mustarstep}) is true. 

$\bullet$ The second equality in (\ref{mustarstep}) can be deduced from the first, by making $n=1$ and $\omega_1=0$.

\begin{subsection}{Bernoulli convolution in base $\beta={1+\sqrt5\over2}$ (\cite{Erd})}
The Gibbs properties of $\mu$ have been studied in \cite{OST3} in the following sense: let be the words
\begin{equation}\label{words}
w(0):=00,\quad w(1)=010\quad\hbox{and}\quad w(2)=10;
\end{equation}
then for any $x\in[0,1[$,
there exists a unique sequence $\xi(x)=(\xi_n)_{n\ge1}\in~\Omega_3$ such that the Parry expansion $\varepsilon(x)$
belongs for any $n\ge1$ to the cylinder $[w(\xi_1\dots\xi_n)]$, where $w(\xi_1\dots\xi_n)$ 
is the concatenation of the words $w(\xi_1),\dots,w(\xi_n)$. The measure
$\mu\circ\xi^{-1}$ is weak Gibbs on
$\Omega_3$ if and only if $p=q$ (this case is studied more in details in \cite{FO}); nevertheless
$\phi_{\mu\circ\xi^{-1}}(100\dots)=\infty$ in this case.
\end{subsection}

\begin{subsection}{Bernoulli convolution in base $-\beta=-{1+\sqrt5\over2}$}
The measure $\mu_\star$ has better Gibbs properties than $\mu$: let us consider now -- for any $x\in[0,1[$ -- the sequence
$\xi_\star(x)=(\xi^\star_n)_{n\ge1}\in\Omega_3$ such that $\alpha(x)\in[w(\xi^\star_1\dots
\xi^\star_n)]_\star$ for all $n\ge1$, we have the following

\begin{theorem}\label{Gibbs*}(i) If $p\ge q$ the measure $\mu_\star\circ\xi_\star^{-1}$ is weak Gibbs on
$\Omega_3$~. 

(ii) if $p\le q$ the measure $\mu_\star\circ S\circ\xi_\star^{-1}$ is weak Gibbs on $\Omega_3$~, where
$S(x)=1-x$ for any $x\in[0,1]$.
\end{theorem}

\begin{proof}
$(ii)$ can be deduced from $(i)$ by using the symmetry relation
$$
Y(\omega_1\omega_2\dots)=1-Y((1-\omega_1)(1-\omega_2)\dots),
$$
which implies $\mu_\star^{(p,q)}\circ S=\mu_\star^{(q,p)}$.

In order to prove $(i)$, we don't use the matrices $A_k$ but the product matrices associated to the three words defined in (\ref{words}): setting
$\alpha={p\over q}$ we~have
$$
A^*_0:={A_0}^2=pq\left(\begin{array}{ccc}\alpha&1&{1\over\alpha}\cr0&1&0\cr0&0&1\end{array}\right),\quad
A^*_1:=A_0A_1A_0=pq^2\left(\begin{array}{ccc}\alpha&1&0\cr\alpha&1&{1\over\alpha}\cr0&0&0\end{array}\right),$$$$
A^*_2:=A_1A_0=pq\left(\begin{array}{ccc}1&{1\over\alpha}&0\cr0&0&0\cr\alpha&1&{1\over\alpha}\end{array}\right).
$$
Let us prove $(i)$ by means of Proposition \ref{pot=>weak}: more precisely we shall prove the uniform convergence of the
(continuous) $n$-step potential $\phi_n:\Omega_3\to{\mathbb R}$ defined by
\begin{equation}\label{phine}
\phi_n(\omega):=\log{\mu_\star\circ\xi_\star^{-1}[\omega_1\dots\omega_n]\over\mu_\star\circ\xi_\star^{-1}[\omega_2\dots\omega_n]}
=\log{\left(\begin{array}{ccc}1&0&0\end{array}\right)A^*_{\omega_1}\dots
A^*_{\omega_n}\left(\begin{array}{c}1
\cr{q\over
1+q}\cr{1\over
1+q}\end{array}\right)\over\left(\begin{array}{ccc}1&0&0\end{array}\right)A^*_{\omega_2}\dots A^*_{\omega_n}\left(\begin{array}{c}1
\cr{q\over
1+q}\cr{1\over
1+q}\end{array}\right)}.
\end{equation}
Notice that
\begin{equation}\label{sanspnniqn}
{A^*_0}^n=(pq)^n\left(\begin{array}{ccc}v_n(\alpha)&\alpha
u_n(\alpha)&u_n(\alpha)\cr0&1&0\cr0&0&1\end{array}\right)\quad\hbox{ and }\quad{A^*_2}^n=(pq)^n\left(\begin{array}{ccc}1&1/\alpha&0\cr0&0&0\cr
u_n(1/\alpha)&u_n(1/\alpha)/\alpha&v_n(1/\alpha)\end{array}\right)
\end{equation}
where $u_n(x):=x^{-1}+x^0+\dots+x^{n-2}$ and $v_n(x):=x^n$ for any positive real~$x$.

\qquad From now on we use the formalism of continued
fractions (\cite{Per}) in a same way as in \cite{OST3}: given $n$ (odd) and
$a_0\ge0$,
$a_1>0,\dots,a_n>0$ we put
$$
\left(\begin{array}{c}p_{-1}\cr q_{-1}\end{array}\right)=\left(\begin{array}{cc}1\cr0\end{array}\right),\quad
\left(\begin{array}{c}p_0\cr q_0\end{array}\right)=\left(\begin{array}{c}u_0\cr1\end{array}\right)\quad\hbox{and, for }1\le k\le n,\ 
\left(\begin{array}{c}p_k\cr q_k\end{array}\right)=u_k\left(\begin{array}{c}p_{k-1}\cr q_{k-1}\end{array}\right)+v_{k-1}\left(\begin{array}{c}p_{k-2}\cr q_{k-2}\end{array}\right)
$$
where, for our purpose, $\left\{\begin{array}{ll}u_i:=u_{a_i}(\alpha)&\hbox{($i$ even)}\cr
u_i:=u_{a_i}(1/\alpha)&\hbox{($i$ odd)}\end{array}\right.$
\quad and\quad$\left\{\begin{array}{ll}v_i:=v_{a_i}(\alpha)&\hbox{($i$ even)}
\cr v_i:=v_{a_i}(1/\alpha)&\hbox{($i$ odd)}.\end{array}\right.$
We have
\begin{equation}\label{avecpnetqn}
{A^*_0}^{a_0}{A^*_2}^{a_1}{A^*_0}^{a_2}\dots{A^*_2}^{a_n}=(pq)^{a_0+\dots+a_n}\left(\begin{array}{ccc}p_n&p_n/\alpha&v_np_{n-1}\cr0&0&0\cr q_n&q_n/\alpha&v_nq_{n-1}\end{array}\right).
\end{equation}

The difference
$\displaystyle\delta_k=~\left\vert{p_k\over q_k}-{p_{k-1}\over q_{k-1}}\right\vert$ is known to be at most
$\displaystyle{1\over a_1+\dots+a_k}$ in the case of the regular continued fractions (\cite{Khin}) that is -- with our notations
-- in the case
$\alpha=1$. We complete by the following
%%%%%%%%%%%%%%%%%%%%%%%
\begin{lemma}\label{annexe}
{\sl  If $\alpha>1$, then

(i) for $k\ge1$,\enskip $\displaystyle\delta_k
\le{v_{k-1}\over u_ku_{k-1}+v_{k-1}}\ \delta_{k-1}$~;

(ii) for $k\ge1$ even,\enskip $\displaystyle\delta_k
\le\alpha^{1-(a_{k-1}+a_k)}\delta_{k-1}$~;

(iii) for $k\ge1$ even,\enskip $\displaystyle\delta_k
\le\alpha^{a_0-(a_1+\dots+a_k)/2}$.}
\end{lemma}
%%%%%%%%%%%%%%%%

\begin{proof} $(i)$ By the definition of $p_k$ and $q_k$~,
%%%%%%%%%%%%%%%%%
\begin{eqnarray*}
{p_k\over q_k}
\!\!\!&=&\!\!\!{u_kp_{k-1}+v_{k-1}p_{k-2}\over
u_kq_{k-1}+v_{k-1}q_{k-2}}\\
\!\!\!&=&\!\!\!{u_kq_{k-1}\over u_kq_{k-1}+v_{k-1}q_{k-2}}\cdot
{p_{k-1}\over
q_{k-1}}+{v_{k-1}q_{k-2}\over u_kq_{k-1}+
v_{k-1}q_{k-2}}\cdot{p_{k-2}\over q_{k-2}}
\end{eqnarray*}
%%%%%%%%%%%%%%%%%%%%%%%
hence
%%%%%%%%%%%%%%%%%
\begin{eqnarray*}
{p_k\over q_k}-{p_{k-1}\over q_{k-1}}
\!\!\!&=&\!\!\!{{v_{k-1}q_{k-2}}\over u_kq_{k-1}+
{v_{k-1}q_{k-2}}}\cdot\left({p_{k-2}\over q_{k-2}}-{p_{k-1}\over q_{k-1}}\right)
\end{eqnarray*}
%%%%%%%%%%%%%%%%%%%%%%%
and, since $q_{k-1}\ge u_{k-1}q_{k-2}$~, we are done.

$(ii)$ If $k$ is even one has $\displaystyle{v_{k-1}\over u_ku_{k-1}}\le{\alpha^{-a_{k-1}}\over\alpha^{a_k-2}\alpha}$~,
hence $(i)$ implies $(ii)$.

$(iii)$ The inequalities $(i)$ and $(ii)$ imply respectively that the sequence $(\delta_k)$ is non-increasing
and, if $k$ is even, $\displaystyle\delta_k
\le\alpha^{-(a_{k-1}+a_k)/2}\delta_{k-1}$~; hence $\displaystyle\delta_2
\le\alpha^{-(a_1+a_2)/2}\delta_1$ and, by induction $\displaystyle\delta_k
\le~\alpha^{-(a_1+\dots+a_k)/2}\delta_1$ for any $k$ even. Now $\displaystyle\delta_1={v_0\over
u_1}\le\alpha^{a_0}$.

\end{proof}

Notice that this lemma implies $\displaystyle\delta_k
\le\alpha^{a_0-(k-1)/2}$ for any $k\ge1$, hence the sequence $\displaystyle k\mapsto{p_k\over q_k}$
converges. Now we can prove the following
%%%%%%%%%%%%%%%%%%%%%%%
\begin{lemma}\label{zerodeux}Suppose $\alpha\ge1$ and let $\omega\in\Omega_3$~.

(i) At least one of the followings assertions is true:

$\exists N\ge0$ such that
$\omega_{N+1}\dots\omega_{N+n}\in\{0,2\}^{n-1}\times\{2\}$ for infinitely many $n\ge1$;

$\exists N\ge0$ such that
$\omega_{N+1}\dots\omega_{N+n}\in\{0\}^n$ for all $n\ge1$;

$\exists N\ge0$ and $n\ge2$ such that
$\omega_{N+1}\dots\omega_{N+n}\in\{1,2\}\times\{0\}^{n-2}\times\{1\}$.

(ii) In all cases there exists $N\ge0$ and $n\ge1$ such that
$$
h,h'\ge N+n,\ 
\omega'\in[\omega_1\dots\omega_{N+n}]\Rightarrow\left\vert\phi_{h'}(\omega')-\phi_h(\omega)\right\vert\le\varepsilon.
$$
\end{lemma}
%%%%%%%%%%%%%%%%

\begin{proof} $(i)$ If there exists $N\ge0$ such that $\sigma^N\omega\in\{0,2\}^{\mathbb N}$, we are in the two first cases. If not,
the digit $1$ occurs infinitely many times in the sequence $\omega$. The second occurrence of $1$ is necessarily preceded
by a word in $\{1,2\}\times\{0\}^k$ for some $k\ge0$, hence we are in the third case.

$(ii)$ Let $N$ and $n$ be as in $(i)$. From (\ref{phine}), for any $h\ge N+n$ and $\omega'\in[\omega_1\dots\omega_{N+n}]$ there exists some
reals $a,b,c,a',b',c',x,y,z$ such that
\begin{equation}\label{ratio}
\phi_h(\omega')=\log{\left(\begin{array}{ccc}a&b&c\end{array}\right)
A^*_{\omega_{N+1}}\dots A^*_{\omega_{N+n}}\left(\begin{array}{c}x\cr y\cr z\end{array}\right)\over
\left(\begin{array}{ccc}a'&b'&c'\end{array}\right)A^*_{\omega_{N+1}}\dots A^*_{\omega_{N+n}}\left(\begin{array}{c}x\cr y\cr z\end{array}\right)},
\end{equation}
where only $x,y$ and $z$ depend on $h$ and $\omega'$.

Suppose the first assertion in $(i)$ is true. We deduce from the expression of $A^*_{\omega_{N+1}}\dots A^*_{\omega_{N+n}}$ in
(\ref{avecpnetqn}) that
$$
\phi_h(\omega')=\log{(ap_n+cq_n)(x+{y\over\alpha})+v_n(ap_{n-1}+cq_{n-1})z\over(a'p_n+c'q_n)(x+{y
\over\alpha})+v_n(a'p_{n-1}+c'q_{n-1})z}.
$$
This ratio lies between $\log\displaystyle{ap_n+cq_n\over a'p_n+c'q_n}$ and $\log\displaystyle{ap_{n-1}+cq_{n-1}\over
a'p_{n-1}+c'q_{n-1}}$. These bounds do not depend on $h$ nor $\omega'$, and converge -- for $n\to\infty$ -- to
$\log\displaystyle{a\theta+c\over a'\theta+c'}$ when
$n\to\infty$, where
$\displaystyle\theta:=\lim_{n\to\infty}{p_n\over q_n}$. We deduce that $(ii)$ is true in this case, by choosing $n$ large enough.

The proof is similar when the second assertion in $(i)$ is true, by using the expression of $A^*_{\omega_{N+1}}\dots A^*_{\omega_{N+n}}$ in
(\ref{sanspnniqn}).

If the third assertion in $(i)$ is true, the matrix $A^*_{\omega_{N+1}}\dots A^*_{\omega_{N+n}}$ has rank $1$; whence it maps ${\mathbb R}^3$ into a
space of dimension $1$, and the ratio in (\ref{ratio}) do not depend on $h$ nor $\omega'$ so that
$$
h,h'\ge N+n,\ 
\omega'\in[\omega_1\dots\omega_{N+n}]\Rightarrow\left\vert\phi_{h'}(\omega')-\phi_h(\omega)\right\vert=0.
$$

\end{proof}

{\bf  End of the proof of Theorem \ref{Gibbs*}. } Notice that Lemma \ref{zerodeux}$(ii)$ implies -- by making $\omega=\omega'$~-- that the sequence
$h\mapsto\phi_h(\omega)$ is Cauchy; let $\phi(\omega)$ be its limit.

Now we make $h,h'\to\infty$ in Lemma \ref{zerodeux}$(ii)$: we obtain
$\left\vert\phi(\omega')-\phi(\omega)\right\vert\le~\varepsilon$ for any $\omega'$ in the neighborhood $[\omega_1\dots\omega_{N+n}]$ of $\omega$, and this
prove the continuity of $\phi$ so, by Proposition \ref{pot=>weak} $\mu_\star\circ\xi_\star^{-1}$ is weak Gibbs.

\end{proof}

\end{subsection}

%%%%%%%%%%%%%%%%%%%%%%%%%%%%%%%%%%%%%%%%
%%%%%%%%%%%%%% BIBLIOGRAPHY %%%%%%%%%%%%%%%%
%%%%%%%%%%%%%%%%%%%%%%%%%%%%%%%%%%%%%%%%


\begin{thebibliography}{99}
%%%%%%%%%%%%%%%%%%%%%%%%%%%%%%%%%%%%%%%%%

%%%%%%%%%%%%%%%%%%%%%%%%%%%%%%%%%%%%%%%%%
\bibitem{Bow1} \textsc{R. Bowen},
\textit{Equilibrium states and the ergodic theory of Anosov diffeomorphisms}.
Lecture Notes in Mathematics {\bf 470}, Berlin-New York, 1975, i+108 pp.
%%%%%%%%%%%%%%%%%%%%%%%%%%%%%%%%%%%%%%%%%
\bibitem{Bow} \textsc{R. Bowen},
\textit{Smooth partitions of Anosov diffeomorphisms are weak Bernoulli}.
Israel J. Math. {\bf 21} (1975), 95--100.
%%%%%%%%%%%%%%%%%%%%%%%%%%%%%%%%%%%%%%%%%
\bibitem{BMP} \textsc{G. Brown, G. Michon and J. Peyri\`ere},
\textit{On the multifractal analysis of measures}.
J. Stat. Phys. {\bf 66} (1992), 775--790.
%%%%%%%%%%%%%%%%%%%%%%%%%%%%%%%%%%%%%%%%%
\bibitem{DS} \textsc{Y. Dupain and V.T. S\'os},
\textit{On the one-sided boundedness of the discrepancy-function of the
sequence $\{n\alpha\}$}.
Acta Arith. {\bf37} (1980), 363--374.
%%%%%%%%%%%%%%%%%%%%%%%%%%%%%%%%%%%%%%%%%
\bibitem{Erd} \textsc{P. Erd\H os},
\textit{On a family of symmetric Bernoulli convolutions}.
Amer. J. Math. {\bf61} (1939), 974--976.
%%%%%%%%%%%%%%%%%%%%%%%%%%%%%%%%%%%%%%%%%
\bibitem{FO} \textsc{D-J. Feng \& E. Olivier},
\textit{Multifractal analysis of weak Gibbs measures and phase 
transition -- application to some Bernoulli convolutions}.
Erg. Th. \& Dyn. Systems {\bf 23} (2003), 1751--1784.
%%%%%%%%%%%%%%%%%%%%%%%%%%%%%%%%%%%%%%%%%
\bibitem{Heu} \textsc{Y. Heurteaux},
\textit{Estimations de la dimension inf\'erieure
et de la dimension sup\'erieure des mesures}.
Ann. Inst. Henri Poincar\'e {\bf 34} (1998), 309--338.
%%%%%%%%%%%%%%%%%%%%%%%%%%%%%%%%%%%%%%%%%
\bibitem{Kea}  \textsc{M. Keane},
\textit{Strongly Mixing $g$-Measures}.
Invent. Math. {\bf 16} (1972), 309--324.
%%%%%%%%%%%%%%%%%%%%%%%%%%%%%%%%%%%%%%%%%
\bibitem{Khin} \textsc{A.Y. Khinchin},
\textit{Continued Fractions}.
Chicago and London, 1964, 95 pp.
%%%%%%%%%%%%%%%%%%%%%%%%%%%%%%%%%%%%%%%%%
\bibitem{MN} \textsc{A. Mukherjea \& A. Nakassis},
\textit{On the continuous singularity of the limit distribution of products of i. i. d. $d\times d$ 
stochastic matrices}.
J. Theor. Probab. {\bf 15} (2002), 903--918.
%%%%%%%%%%%%%%%%%%%%%%%%%%%%%%%%%%%%%%%%%
\bibitem{MNR} \textsc{A. Mukherjea, A. Nakassis \& J.S. Ratti},
\textit{On the distribution of the limit of products of i. i. d. $2\times2$ random stochastic matrices}. 
J. Theor. Probab. {\bf 12} (1999), 571--583.
%%%%%%%%%%%%%%%%%%%%%%%%%%%%%%%%%%%%%%%%%
\bibitem{OST2} \textsc{E. Olivier, N. Sidorov, \& A. Thomas},
\textit{On the Gibbs properties of Bernoulli convolutions 
and related problems in Fractal Geometry}.
Preprint LATP 02-14 (2002).
%%%%%%%%%%%%%%%%%%%%%%%%%%%%%%%%%%%%%%%%%
\bibitem{OST3} \textsc{E. Olivier, N. Sidorov, \& A. Thomas},
\textit{On the Gibbs properties of Bernoulli convolutions 
related to $\beta$-numeration in multinacci bases}.
Monatshefte f\"ur Math. {\bf 145} (2005), 145--174.
%%%%%%%%%%%%%%%%%%%%%%%%%%%%%%%%%%%%%%%%%
\bibitem{OST} \textsc{A. Ostrowski},
\textit{Bemerkungen zur Theorie der Diophantischen Approximationen I, II}.
Abh. Math. Sem. Hamburg I
(1922), 77--98 and 250--251. 
%%%%%%%%%%%%%%%%%%%%%%%%%%%%%%%%%%%%%%%%%
\bibitem{P-P-W} \textsc{M.R. Palmer, W. Parry \& P. Walters},
\textit{Large Sets of Endomorphisms and of $g$-Measures}.
Lecture Notes in Math. {\bf 668} (1978) 191--210.
%%%%%%%%%%%%%%%%%%%%%%%%%%%%%%%%%%%%%%%%%
\bibitem{P} \textsc{W. Parry},
\textit{On the $\beta$-expansions of real numbers}.
Acta Math. Acad. Sci. Hung. {\bf11} (1960), 401--416.
%%%%%%%%%%%%%%%%%%%%%%%%%%%%%%%%%%%%%%%%%
\bibitem{PS} \textsc{Y. Peres, W. Schlag, B. Solomyak},
\textit{Sixty years of Bernoulli convolutions}.
Prog. Probab. {\bf46} (2000), 39--65.
%%%%%%%%%%%%%%%%%%%%%%%%%%%%%%%%%%%%%%%%%
\bibitem{Per} \textsc{O. Perron},
\textit{Die Lehre von den Kettenbr\"uchen}.
New York, 1950, 524 pp.
%%%%%%%%%%%%%%%%%%%%%%%%%%%%%%%%%%%%%%%%%
\bibitem{SV} \textsc{N. Sidorov, A. Vershik}, 
\textit{Ergodic Properties of the Erd\H os measure, the Entropy
of the Goldenshift, and Related Problems}. 
Monatshefte f\"ur Math. {\bf 126} (1998), 215--261.
%%%%%%%%%%%%%%%%%%%%%%%%%%%%%%%%%%%%%%%%%
\bibitem{Wal} \textsc{P. Walters}, 
\textit{Ruelle's operator theorem and $g$-measures}. 
Trans. Am. Math. Soc. {\bf214} (1975), 375--387. 
%%%%%%%%%%%%%%%%%%%%%%%%%%%%%%%%%%%%%%%%%
\bibitem{Yur} \textsc{M. Yuri},
\textit{Zeta functions for certain non-hyperbolic systems and topological Markov
approximations}.
Erg. Th. \& Dyn. Syst. {\bf18} (1998), 1589--1612.
%%%%%%%%%%%%%%%%%%%%%%%%%%%%%%%%%%%%%%%%%

\end{thebibliography}
\end{document}